\documentclass[11pt,a4paper]{article}
\usepackage[margin=1in]{geometry}
\usepackage[utf8]{inputenc}

\usepackage{tikz}
\usetikzlibrary{positioning,math}

\usepackage{graphicx}

\usepackage{latexsym, amsfonts, amsmath, amssymb, graphicx}
\newcommand{\be}{\begin{equation}}
\newcommand{\ee}{\end{equation}}
\newcommand{\bea}{\begin{eqnarray}}
\newcommand{\eea}{\end{eqnarray}}
\newcommand{\bean}{\begin{eqnarray*}} 
\newcommand{\eean}{\end{eqnarray*}}

\newcommand{\brray}{\begin{array}}
\newcommand{\erray}{\end{array}}
\newcommand{\ben}{\begin{equation}{nonumber}}
\newcommand{\een}{\end{equation}{nonumber}}

\newtheorem{dfn}{Definition}[section]
\newtheorem{thm}[dfn]{Theorem}
\newtheorem{lmma}[dfn]{Lemma}
\newtheorem{ppsn}[dfn]{Proposition}
\newtheorem{crlre}[dfn]{Corollary}
\newtheorem{xmpl}[dfn]{Example}
\newtheorem{rmrk}[dfn]{Remark}
\newcommand{\bdfn}{\begin{dfn}}
\newcommand{\bthm}{\begin{thm}}
\newcommand{\blmma}{\begin{lmma}}
\newcommand{\bppsn}{\begin{ppsn}}
\newcommand{\bcrlre}{\begin{crlre}}
\newcommand{\bxmpl}{\begin{xmpl}}
\newcommand{\brmrk}{\begin{rmrk}}
\newcommand{\edfn}{\end{dfn}}
\newcommand{\ethm}{\end{thm}}
\newcommand{\elmma}{\end{lmma}}
\newcommand{\eppsn}{\end{ppsn}}
\newcommand{\ecrlre}{\end{crlre}}
\newcommand{\exmpl}{\end{xmpl}}
\newcommand{\ermrk}{\end{rmrk}}
\newcommand{\IC}{\mathbb{C}}
\newcommand{\IR}{\mathbb{R}}

\newcommand{\IT}{\mathbb{T}}

\newcommand{\IZ}{\mathbb{Z}}
\newcommand{\cla}{{\cal A}}
\newcommand{\clb}{{\cal B}}
\newcommand{\clc}{{\cal C}}

\newcommand{\cle}{{\cal E}}
\newcommand{\clf}{{\cal F}}
\newcommand{\clg}{{\cal G}}
\newcommand{\clh}{{\cal H}}

\newcommand{\clj}{{\cal J}}

\newcommand{\cln}{{\cal N}}

\newcommand{\clq}{{\cal Q}}

\newcommand{\clt}{{\cal T}}

\newcommand{\clz}{{\cal Z}}

\def\a*{{\cal A}_{h,*}}
\def\B{{\cal B}(h)}
\def\B1{{\cal B}_1(h)}
\def\b{{\cal B}^{\rm s.a.}(h)}
\def\b1{{\cal B}^{\rm s.a.}_1(h)}

\newcommand{\ot}{\otimes}

\newcommand{\raro}{\rightarrow}

\def \qed {$\Box$}

\def\a*{{\cal A}_{h,*}}
\def\B{{\cal B}(h)}
\def\B1{{\cal B}_1(h)}
\def\b{{\cal B}^{\rm s.a.}(h)}
\def\b1{{\cal B}^{\rm s.a.}_1(h)}

\newcommand{\RNum}[1]{\uppercase\expandafter{\romannumeral #1\relax}}

\date{}
\begin{document}

\title{Invariance of KMS states on graph $C^{\ast}$-algebras under classical and quantum symmetry}
\author{Soumalya Joardar, Arnab Mandal}
\maketitle

%\vspace{0.1 in}

\begin{abstract}
 We study invariance of KMS states on graph $C^{\ast}$-algebras coming from strongly connected and circulant graphs under the classical and quantum symmetry of the graphs. We show that  the unique KMS state for strongly connected graphs is invariant under quantum automorphism group of the graph. For circulant graphs, it is shown that  the action of classical and quantum automorphism group preserves only one of the  KMS states occurring at the critical inverse temperature. We also give an example of a graph $C^{\ast}$-algebra having more than one KMS state such that  all of them are invariant under the action of classical automorphism group of the graph, but there is a unique KMS state which is invariant under the action of quantum automorphism group of the graph. 
\end{abstract}
\section {Introduction}
For a $C^{\ast}$-algebra $\cla$, $(\cla,\sigma)$ is called a $C^{\ast}$-dynamical system if there is a strongly continuous map $\sigma:\mathbb{R}\raro{\rm Aut}(\cla)$. When $\clh$ is a finite dimensional Hilbert space, a $C^{\ast}$-dynamical system $(\clb(\clh),\sigma)$ is given by a self adjoint operator $H\in\clb(\clh)$ in the sense that $\sigma_{t}(A)=e^{itH}Ae^{-itH}$. For such a $C^{\ast}$-dynamical system, it is well known that at any inverse temperature $\beta\in\mathbb{R}$, the unique thermal equilibrium state is given by the Gibbs state 
\begin{displaymath}
\omega_{\beta}(A)=\frac{{\rm Tr}(e^{-\beta H}A)}{{\rm Tr}(e^{-\beta H})}, A\in\cla.
\end{displaymath} 
For a general $C^{\ast}$-dynamical system $(\cla,\sigma)$ the generalization of the Gibbs states are the KMS (Kubo-Martin-Schwinger) states. A KMS state for a $C^{\ast}$-dynamical system $(\cla,\sigma)$ at an inverse temperature $\beta\in\mathbb{R}$ is a state $\tau\in\cla^{\ast}$ that satisfies the KMS condition given by
\begin{displaymath}
\tau(ab)=\tau(b\sigma_{i\beta}(a)),
\end{displaymath}
for $a,b$ in a dense subalgebra of $\cla$ called the algebra of analytic elements of $(\cla,\sigma)$. For a general $C^{\ast}$-dynamical system, unlike the finite dimensional case, KMS state does not exist at every temperature. Even though they exist, generally nothing can be said about their uniqueness at a given temperature. It is worth mentioning that in Physics literature, uniqueness of KMS state often is related to phase transition and symmetry breaking.\\
\indent One of the mathematically well studied KMS states are the KMS states for the $C^{\ast}$-dynamical systems on the graph $C^{\ast}$-algebras (see \cite{Laca}, \cite{watatani}). For a finite directed graph $\Gamma$ the dynamical system is given by the $C^{\ast}$-algebra $C^{\ast}(\Gamma)$ and the automorphism $\sigma$ being the natural lift of the canonical gauge action on $C^{\ast}(\Gamma)$. In \cite{Laca} it is shown that there is a KMS state at the inverse critical temperature ${\rm ln}(\rho(D))$ if and only if $\rho(D)$ is an eigen value of $D$ with eigen vectors having all entries non negative, where $\rho(D)$ is the spectral radius of the vertex matrix $D$ of the graph. For a general graph $C^{\ast}$-algebra, we can not say anything about uniqueness of KMS states.  In \cite{Laca}, for strongly connected graphs, such uniqueness result has been obtained. In fact, for strongly connected graphs, there is a unique KMS state occurring only at the critical inverse temperature. However, for another class called the circulant graphs, one can show that KMS states at the critical inverse temperature are not unique. In this context, it is interesting to study invariance of KMS states under some natural added (apart from the gauge symmetry) internal symmetry of the graph $C^{\ast}$-algebra and see if such invariance could force the KMS states to be unique in certain cases. It is shown in \cite{Web} that for a graph $\Gamma$, the graph $C^{\ast}$-algebra has a natural generalized symmetry coming from the quantum automorphism group $Q^{\rm aut}_{\rm Ban}(\Gamma)$ (see \cite{Ban}) of the graph itself. This symmetry is generalized in the sense that it contains the classical automorphism group ${\rm Aut}(\Gamma)$ of the graph. In this paper we study the invariance of the KMS states under this generalized symmetry.\\
\indent For a strongly connected graph $\Gamma$, we show that the already unique KMS state is preserved by the quantum automorphism group $Q_{\rm Ban}^{\rm aut}(\Gamma)$. This result has a rather interesting consequence on the ergodicity of the action of $Q^{\rm aut}_{\rm Ban}(\Gamma)$ on the graph. It is shown that if the row sums of the vertex matrix of a strongly connected graph are not equal, then $Q^{\rm aut}_{\rm Ban}(\Gamma)$ can not act ergodically. Then we study another class of graphs called the circulant graphs. Circulant graphs admit KMS states at the inverse critical temperature, but they are not necessarily unique. But due to the transitivity of the action of the automorphism group, it is shown that in fact there exists a unique KMS state which is invariant under the classical or quantum symmetry of the system. In fact we also show that the only temperature where the KMS state could occur is the inverse critical temperature. Finally, we show by an example that we necessarily need the invariance of the KMS state under the action of the quantum symmetry group to force the KMS state to be unique. More precisely, we give an example of a graph with 48 vertices coming from the Linear Binary Constraint system (LBCS, see \cite{nonlocal1}) where the corresponding graph $C^{\ast}$-algebra has more than one KMS state all of which are preserved by the action of the classical automorphism group of the graph. However, it has a unique $Q^{\rm aut}_{\rm Ban}(\Gamma)$ invariant KMS state. In this example also, the only possible inverse temperature at which the KMS state could occur is the inverse critical temperature. This shows that in deed where the classical symmetry fails to fix KMS state, the richer `genuine' quantum symmetry of the system plays a crucial role to fix KMS state.
\section{Preliminaries}
\subsection{KMS states on graph $C^{\ast}$-algebra without sink}
\label{KMSgraph}
A {\bf finite} directed graph is a collection of finitely many edges and vertices. If we denote the edge set of a graph $\Gamma$ by $E=(e_{1},\ldots,e_{n})$ and set of vertices of $\Gamma$ by $V=(v_{1},\ldots,v_{m})$ then recall the maps $s,t:E\raro V$ and the vertex matrix $D$ which is an $m\times m$ matrix whose $ij$-th entry is $k$ if there are $k$-number of edges from $v_{i}$ to $v_{j}$. We denote the space of paths by $E^{\ast}$ (see \cite{Laca}). $vE^{\ast}w$ will denote the set of paths between two vertices $v$ and $w$.
\bdfn
$\Gamma$ is said to be without sink if the map $s:E\raro V$ is surjective. Furthermore $\Gamma$ is said to be without multiple edges if the adjacency matrix $D$ has entries either $1$ or $0$. 
\edfn
\brmrk
Note that the graph $C^{\ast}$-algebra corresponding to a graph without sink is a Cuntz-Krieger algebra. Reader might see \cite{Cuntz-Krieger} for more details on Cuntz-Krieger algebra.
\ermrk
Now we recall some basic facts about graph $C^{\ast}$-algebras. Reader might consult \cite{Pask} for details on graph $C^{\ast}$-algebras. Let $\Gamma=\{E=(e_{1},...,e_{n}),V=(v_{1},...,v_{m})\}$ be a finite, directed graph without sink. In this paper all the graphs are {\bf finite}, without {\bf sink} and without {\bf multiple edges}. We assign partial isometries $S_{i}$'s to edges $e_{i}$ for all $i=1,...,n$ and projections $p_{v_{i}}$ to the vertices $v_{i}$ for all $i=1,...,m$.
\bdfn
\label{Graph}
The graph $C^{\ast}$-algebra $C^{\ast}(\Gamma)$ is defined as the universal $C^{\ast}$-algebra generated by partial isometries $\{S_{i}\}_{i=1,\ldots,n}$ and mutually orthogonal projections $\{p_{v_{k}}\}_{k=1,\ldots,m}$ satisfying the following relations:
\begin{displaymath}
  S_{i}^{\ast}S_{i}=p_{t(e_{i})}, \sum_{s(e_{j})=v_{l}}S_{j}S_{j}^{\ast}=p_{v_{l}}.
	\end{displaymath}
\edfn
In a graph $C^{\ast}$-algebra $C^{\ast}(\Gamma)$, we have the following (see Subsection 2.1 of \cite{Pask}):\\
1. $\sum_{k=1}^{m}p_{v_{k}}=1$ and $S_{i}^{\ast}S_{j}=0$ for $i\neq j$.
\vspace{0.05in}\\
2. $S_{\mu}=S_{1}S_{2}\ldots S_{l}$ is non zero if and only if $\mu=e_{1}e_{2}\ldots e_{l}$ is a path i.e. $t(e_{i})=s(e_{i+1})$ for $i=1,\ldots,(l-1)$.
\vspace{0.05in}\\
3. $C^{\ast}(\Gamma)={\overline{\rm Sp}}\{S_{\mu}S_{\nu}^{\ast}:t(\mu)=t(\nu)\}$. \vspace{0.1in}\\
\indent Now we shall briefly discuss KMS-states on graph $C^{\ast}$-algebras coming from graphs without sink. For that we recall Toeplitz algebra $\clt C^{\ast}(\Gamma)$. Readers are referred to \cite{Laca} for details. Our convention though is opposite to that of \cite{Laca} in the sense that we interchange source projections and target projections. Also we shall modify the results of \cite{Laca} according to our need. Suppose that $\Gamma$ is a directed graph as before. A Toeplitz-Cuntz-Krieger $\Gamma$ family consists of mutually orthogonal projections $\{p_{v_{i}}:v_{i}\in V\}$ and partial isometries $\{S_{i}:e_{i}\in E\}$ such that $\{S_{i}^{\ast}S_{i}=p_{t(e_{i})}\}$ and
\begin{displaymath}
p_{v_{l}}\geq \sum_{s(e_{i})=v_{l}}S_{i}S_{i}^{\ast}.
\end{displaymath}
Toeplitz algebra $\clt C^{\ast}(\Gamma)$ is defined to be the universal $C^{\ast}$-algebra generated by the Toeplitz-Cuntz-Krieger $\Gamma$ family. It is clear from the definition that $C^{\ast}(\Gamma)$ is the quotient of $\clt C^{\ast}(\Gamma)$ by the ideal $\clj$ generated by
\begin{displaymath}
  P:=\{p_{v_{l}}-\sum_{s(e_{i})=v_{l}}S_{i}S_{i}^{\ast}\}.
	\end{displaymath}
	The standard arguments give $\clt C^{\ast}(\Gamma)=\overline{\rm Sp} \{S_{\mu}S_{\nu}^{\ast}:t(\mu)=t(\nu)\}$. $\clt C^{\ast}(\Gamma)$ admits the usual gauge action $\gamma$ of $\mathbb{T}$ which descends to the usual gauge action on $C^{\ast}(\Gamma)$ given on the generators by $\gamma_{z}(S_{\mu}S_{\nu}^{\ast})=z^{(|\mu|-|\nu|)}S_{\mu}S_{\nu}^{\ast}$. Consequently it has a dynamics $\alpha:\mathbb{R}\raro {\rm Aut} \ C^{\ast}(\Gamma)$ which is lifted from $\gamma$ via the map $t\raro e^{it}$. We recall the following from \cite{Laca} (Proposition 2.1).
	\bppsn
	\label{exist_KMS}
	Let $\Gamma$ be a finite, directed, connected graph without sink and $\gamma:\mathbb{T}\raro {\rm Aut} \ \clt C^{\ast}(\Gamma)$ be the gauge action with the corresponding dynamics $\alpha:\mathbb{R}\raro {\rm Aut} \ \clt C^{\ast}(\Gamma)$. Let $\beta\in\mathbb{R}$.\\
	(a) A state $\tau$ is a ${\rm KMS}_{\beta}$ state of $(\clt C^{\ast}(\Gamma),\alpha)$ if and only if
	\begin{displaymath}
	\tau(S_{\mu}S_{\nu}^{\ast})=\delta_{\mu,\nu}e^{-\beta|\mu|}\tau(p_{t(\mu)}).
	\end{displaymath}
	(b) Suppose that $\tau$ is a ${\rm KMS}_{\beta}$ state of $(\clt C^{\ast}(\Gamma),\alpha)$, and define $\cln^{\tau}=(\cln^{\tau}_{i})\in[0,\infty)^{m}$ by $\cln^{\tau}_{i}=\phi(p_{v_{i}})$. Then $\cln^{\tau}$ is a probability measure on $V$ satisfying the subinvariance condition $D\cln^{\tau}\leq e^{\beta}\cln^{\tau}$.\\
	(c) A ${\rm KMS}_{\beta}$ state factors through $C^{\ast}(\Gamma)$ if and only if $(D\cln^{\tau})_{i}=e^{\beta}\cln^{\tau}_{i}$ for all $i=1,\ldots, m$ i.e. $\cln^{\tau}$ is an eigen vector of $D$ with eigen value $e^{\beta}$.
	\eppsn
\subsubsection{KMS states at critical inverse temperature}
\label{KMS}
In this subsection we collect a few results on existence of KMS states at inverse critical temperature on graph $C^{\ast}$-algebras coming from graphs {\bf without} sink. For that we continue to assume $\Gamma$ to be a finite, connected graph without sink and with vertex matrix $D$. We denote the spectral radius of $D$ by $\rho(D)$. With these notations, Combining Proposition 4.1 and Corollary 4.2 of \cite{Laca}, we have the following
\bppsn
\label{exist1}
The graph $C^{\ast}$-algebra $C^{\ast}(\Gamma)$ has a ${\rm KMS}_{{\rm ln}(\rho(D))}$ state if and only if $\rho(D)$ is an eigen value of $D$ such that it has eigen vector with all entries non negative. 
\eppsn 
\blmma
\label{onetemp}
Suppose $\Gamma$ is a finite directed graph without sink with vertex matrix $D$. If $\rho(D)$ is an eigen value of $D$ with an eigen vector ${\bf v}$ whose entries are strictly positive such that ${\bf v}^{T}D=\rho(D){\bf v}^{T}$, then the only possible inverse temperature where the KMS state could occur is ${\rm ln}(\rho(D))$.
\elmma
{\it Proof}:\\
Suppose $\beta\in\mathbb{R}$ be another possible inverse temperature where a KMS state say $\phi$ could occur. Then since we have assumed our graph to be without sink, by (c) of Proposition \ref{exist_KMS}, $e^{\beta}$ is an eigen value of $D$. Let us denote an eigen vector corresponding to $e^{\beta}$ by ${\bf w}=(w_{1},\ldots,w_{m})$ so that $w_{i}=\phi(p_{v_{i}})$. Since $\phi$ is a state, $w_{i}\geq 0$ for all $i=1,\ldots,m$ with atleast one entry strictly positive. We have
\begin{eqnarray*}
	&& D{\bf w}=e^{\beta}{\bf w}\\
	&\Rightarrow& {\bf v}^{T}D{\bf w}={\bf v}^{T}e^{\beta}{\bf w}\\
	& \Rightarrow&  (\rho(D)-e^{\beta}){\bf v}^{T}{\bf w}=0\\
	\end{eqnarray*}  
By assumption, all the entries of ${\bf v}$ are strictly positive and $w_{i}\geq 0$ for all $i$ with at least one entry strictly positive which imply that ${\bf v}^{T}{\bf w}\neq0$ and hence $e^{\beta}=\rho(D)$ i.e. $\beta={\rm ln}(\rho(D))$.\qed\vspace{0.1in}\\
We discuss examples of two classes of graphs which are {\bf without sink} such that they admit KMS states only at inverse critical temperature. We shall use them later in this paper. 
\bxmpl
\rm{{\bf Strongly connected graphs}: 
\bdfn
A graph is said to be strongly connected if $vE^{\ast}w$ is non empty for all $v,w\in V$.
\edfn
\bdfn
An $m\times m$ matrix $D$ is said to be irreducible if for $i,j\in\{1,\ldots,m\}$, there is some $k>0$ such that $D^{k}(i,j)>0$.
\edfn
\bppsn
A graph is strongly connected if and only if its vertex matrix is irreducible.
\eppsn
\bppsn
 An irreducible matrix $D$ has its spectral radius $\rho(D)$ as an eigen value with one dimensional eigen space spanned by a vector with all its entries strictly positive (called the Perron-Frobenius eigen vector).
\eppsn
As a corollary we have
\bcrlre
Let $\Gamma$ be a strongly connected graph. Then the graph $C^{\ast}$-algebra $C^{\ast}(\Gamma)$ has a unique ${\rm KMS}_{{\rm ln}(\rho(D))}$ state. In fact by (b) of Theorem 4.3 of \cite{Laca}, this is the only KMS state.
\ecrlre}\exmpl
\bxmpl\label{circulant}\rm{{\bf Circulant graphs}
\bdfn
A graph with $m$ vertices is said to be circulant if its automorphism group contains the cyclic group $\mathbb{Z}_{m}.$
\edfn
It is easy to see that if a graph is circulant, then its vertex matrix is determined by its first row vector say $(d_{0},\ldots,d_{m-1})$. More precisely the vertex matrix $D$ of a circulant graph is given by
\begin{center}
$\begin{bmatrix}
d_{0} & d_{1}... & d_{m-1}\\
d_{m-1} & d_{0}... & d_{m-2}\\
.\\
.\\
.\\
d_{1} & d_{2}... & d_{0}
\end{bmatrix}.$
\end{center}
\brmrk
\label{cws}
Note that a circulant graph is always without sink except the trivial case where it has no edge at all. This is because if $i$-th vertex of a circulant matrix is a sink, then the $i$-th row of the vertex matrix will be zero forcing all the rows to be identically zero. For this reason in this paper we study the circulant graphs {\bf without} sink.
\ermrk 
 Let $\epsilon$ be a primitive $m$-th root of unity. The following is well known (see \cite{note}):
\bppsn
For a circulant graph with vertex matrix as above, the eigen values are given by
\begin{displaymath}
\lambda_{l}=d_{0}+\epsilon^{l}d_{1}+\ldots+\epsilon^{(m-1)l}d_{m-1},l=0,\ldots,(m-1)
\end{displaymath}
\eppsn
It is easy to see that $\lambda=\sum_{i=0}^{m-1}d_{i}$ is an eigen value of $D$ and it has a normalized eigen vector given by $(\frac{1}{m},\ldots,\frac{1}{m})$. Since $|\lambda_{l}|\leq \lambda$, we have
\bcrlre
For a circulant graph $\Gamma$ with vertex matrix $D$, $D$ has its spectral radius $\lambda$ as an eigen value with a normalized eigen vector (not necessarily unique) having all its entries non negative. 
\ecrlre
Combining the above corollary with Proposition \ref{exist1}, we have
\bcrlre
\label{not_unique}
For a circulant graph $\Gamma$, $C^{\ast}(\Gamma)$ has a ${\rm KMS}_{{\rm ln}(\lambda)}$ state.
\ecrlre
\blmma
\label{circonetemp}
For a circulant graph $\Gamma$, the only possible temperature where a KMS state could occur is the inverse critical temperature.
\elmma
{\it Proof}:\\
As we have assumed the circulant graphs are without sink, it is enough to show that the vertex matrix $D$ of $\Gamma$ satisfies the conditions of Lemma \ref{onetemp}. It is already observed that the eigen value $\lambda$ has an eigen vector with all its entries positive (column vector with all its entries $1$ to be precise). Also since the row sums are equal to column sums which are equal to $\lambda$, $(1,\ldots,1)D=\lambda(1,\ldots,1)$. Hence an application of Lemma \ref{onetemp} finishes the proof.\qed\vspace{0.1in}\\
Note that  KMS state at inverse critical temperature is not necessarily unique, since the dimension of the eigen space of the eigen value $\lambda$ could be strictly larger than 1 as the following example illustrates:\\
We take the graph whose vertex matrix is given by 
\begin{center}
$\begin{bmatrix}
1 & 0 & 1 & 0\\
0 & 1 & 0 & 1\\
1 & 0 & 1 & 0\\
0 & 1 & 0 & 1
\end{bmatrix}. $
\end{center}
Hence the graph is circulant with its spectral radius $2$ as an eigen value with multiplicity $2$. So the dimension of the corresponding eigen space is $2$ violating the uniqueness of the KMS state at the inverse critical temperature ${\rm ln}(2)$. 
%\begin{figure}[h] 
	%\centering
	%\includegraphics[width=0.4\textwidth]{drawing.pdf}
	%\caption{}
%\end{figure}
%\brmrk
%\label{onetemp}
%It can be shown that the only possible inverse temperature where KMS state occur for circulant graphs is the inverse critical temperature. In deed, for $\beta<{\rm ln}(\lambda)$ using the same argument used to prove (c) of Theorem 4.3 of \cite{Laca}, it can be shown that there is no ${\rm KMS}_{\beta}$ state. Since $\lambda\geq 1$, for $\beta>{\rm ln}(\lambda)$, $e^{\beta}>\lambda$. But $\lambda$ being the spectral radius of the vertex matrix, any $e^{\beta}>\lambda$ can not be an eigen value of $D$. Hence by (c) of Proposition \ref{exist_KMS}, there is no ${\rm KMS}_{\beta}$ state at $\beta>{\rm ln}(\lambda)$.
%\ermrk
}\exmpl

\subsection{Quantum automorphism group of graphs as symmetry of graph $C^{\ast}$-algebra}
\subsubsection{Compact quantum groups and quantum automorphism groups}
\label{qaut}
In this subsection we recall the basics of compact quantum groups and their actions on $C^{\ast}$-algebras. The facts collected in this Subsection are well known and we refer the readers to \cite{Van}, \cite{Woro}, \cite{Wang} for details. All the tensor products in this paper are minimal.

\bdfn
A compact quantum group (CQG) is a pair $(\clq,\Delta)$ such that $\clq$ is a unital $C^{\ast}$-algebra and $\Delta:\clq\raro\clq\ot\clq$ is a unital $C^{\ast}$-homomorphism satisfying\\
(i) $({\rm id}\ot\Delta)\circ\Delta=(\Delta\ot{\rm id})\circ\Delta$.\\
(ii) {\rm Span}$\{\Delta(\clq)(1\ot\clq)\}$ and {\rm Span}$\{\Delta(\clq)(\clq\ot 1)\}$ are dense in $\clq\ot\clq$.
\edfn
Given a CQG $\clq$, there is a canonical dense Hopf $\ast$-algebra $\clq_{0}$ in $\clq$ on which an antipode $\kappa$ and counit $\epsilon$ are defined. Given two CQG's $(\clq_{1},\Delta_{1})$ and $(\clq_{2},\Delta_{2})$, a CQG morphism between them is a $C^{\ast}$-homomorphism $\pi:\clq_{1}\raro\clq_{2}$ such that $(\pi\ot\pi)\circ\Delta_{1}=\Delta_{2}\circ\pi$.
%3. $A_{aut}(M_{n})$ be the universal $C^{\ast}$ algebra generated by $\{q_{ij}^{kl}\}_{i,j,k,l=1,\ldots,n}$ satisfying the following relations:
%\begin{eqnarray*}
%&& \sum_{v=1}^{n}q_{ij}^{kv}q_{rs}^{vl}=\delta_{jr}q_{is}^{kl}, \ i,j,k,l,r,s=1,\ldots,n\\
%&& \sum_{v=1}^{n}q_{lv}^{sr}q_{vk}^{ji}=\delta_{jr}q^{si}_{lk}, \ i,j,k,l,r,s=1,\ldots,n\\
%&& q^{kl^{\ast}}_{ij}=q_{ji}^{lk}, \ i,j,k=1,\ldots,n\\
%&& \sum_{r=1}^{n}q_{rr}^{kl}=\sum_{r=1}^{n}q_{kl}^{rr}=\delta_{kl}, \ k,l=1,\ldots,n.
%\end{eqnarray*}
%Then $A_{aut}(M_{n})$ is a CQG where the coproduct is given by $\Delta(q^{kl}_{ij})=\sum_{r,s=1}^{n}q_{rs}^{kl}\ot q_{ij}^{rs}$ (see Theorem 4.1 of \cite{Wang}).\\

\bdfn
Given a (unital) $C^{\ast}$-algebra $\clc$, a CQG $(\clq,\Delta)$ is said to act faithfully on $\clc$ if there is a unital $C^{\ast}$-homomorphism $\alpha:\clc\raro\clc\ot\clq$ satisfying\\
(i) $(\alpha\ot {\rm id})\circ\alpha=({\rm id}\ot \Delta)\circ\alpha$.\\
(ii) {\rm Span}$\{\alpha(\clc)(1\ot\clq)\}$ is dense in $\clc\ot\clq$.\\
(iii) The $\ast$-algebra generated by the set  $\{(\omega\ot{\rm id})\circ\alpha(\clc): \omega\in\clc^{\ast}\}$ is norm-dense in $\clq$.
\edfn
\bdfn
An action $\alpha:\clc\raro\clc\ot\clq$ is said to be ergodic if $\alpha(c)=c\ot 1$ implies $c\in\mathbb{C}1$.
\edfn
\bdfn
\label{statepres}
Given an action $\alpha$ of a CQG $\clq$ on a $C^{\ast}$-algebra $\clc$, $\alpha$ is said to preserve a state $\tau$ on $\clc$ if $(\tau\ot{\rm id})\circ\alpha(a)=\tau(a)1$ for all $a\in\clc$.
\edfn
For a faithful action of a CQG $(\clq,\Delta)$ on a unital $C^{\ast}$-algebra $\clc$, there is a norm dense $\ast$-subalgebra $\clc_{0}$ of $\clc$ such that the canonical Hopf-algebra $\clq_{0}$ coacts on $\clc_{0}$.
\bdfn
(Def 2.1 of \cite{Bichon})
Given a unital $C^{\ast}$-algebra $\clc$, quantum automorphism group of $\clc$ is a CQG $(\clq,\Delta)$ acting faithfully on $\clc$ satisfying the following universal property:\\
\indent If $(\clb,\Delta_{\clb})$ is any CQG acting faithfully on $\clc$, there is a surjective CQG morphism $\pi:\clq\raro\clb$ such that $({\rm id}\ot \pi)\circ\alpha=\beta$, where $\beta:\clc\raro\clc\ot\clb$ is the corresponding action of $(\clb,\Delta_{\clb})$ on $\clc$ and $\alpha$ is the action of $(\clq,\Delta)$ on $\clc$.
\edfn
\brmrk
In general the universal object might fail to exist in the above category. To ensure existence one generally assumes that the action preserves some fixed state on the $C^{\ast}$-algebra in the sense  of definition \ref{statepres}. We would not go into further details as we are not going to use it in this paper. For further details, reader might consult \cite{Wang}.
\ermrk
\bxmpl
\label{S}
\rm{If we take a space of $n$ points $X_{n}$ then the quantum automorphism group of the $C^{\ast}$-algebra $C(X_{n})$ is given by the CQG (denoted by $S_{n}^{+}$) which as a $C^{\ast}$-algebra is the universal $C^{\ast}$ algebra generated by $\{u_{ij}\}_{i,j=1,\ldots,n}$ satisfying the following relations (see Theorem 3.1 of \cite{Wang}):
	\begin{displaymath}
	u_{ij}^{2}=u_{ij},u_{ij}^{\ast}=u_{ij},\sum_{k=1}^{n}u_{ik}=\sum_{k=1}^{n}u_{ki}=1, i,j=1,\ldots,n.
	\end{displaymath}
	The coproduct on the generators is given by $\Delta(u_{ij})=\sum_{k=1}^{n}u_{ik}\ot u_{kj}$.}\exmpl
%\bxmpl
%\rm{If we take the $C^{\ast}$-algebra $M_{n}(\mathbb{C})$, then as remarked earlier if we take the category of CQG's only having faithful $C^{\ast}$-action on $M_{n}(\mathbb{C})$, then universal object fails to exist in the category. But in addition if we assume any object in the category also has to preserve a linear functional $\phi$ on $M_{n}(\mathbb{C})$, then existence of the universal object can be shown (see \cite{Wang}).}\exmpl 

\subsubsection{Quantum automorphism group of finite graphs and graph $C^{\ast}$-algebras}
Recall the definition of finite, directed graph $\Gamma=((V=v_{1},\ldots,v_{m}),( E= e_{1},\ldots, e_{n}))$ without multiple edge from Subsection \ref{KMSgraph} and the Example \ref{S} of the CQG $S_{n}^{+}$.
\bdfn
\label{qsymban}
$(Q^{\rm aut}_{\rm Ban}(\Gamma),\Delta)$ for a graph $\Gamma$ without multiple edge  is defined to be the quotient $S^{+}_{n}/(AD-DA)$, where $A=((u_{ij}))_{i.j=1,\ldots,m}$, and $D$ is the adjacency matrix for $\Gamma$. The coproduct on the generators is  given by $\Delta(u_{ij})=\sum_{k=1}^{m}u_{ik}\ot u_{kj}$.
\edfn
For the classical automophism group ${\rm Aut}(\Gamma)$, the commutative CQG $C({\rm Aut}(\Gamma))$ is generated by $u_{ij}$ where $u_{ij}$ is a function on $S_{n}$ taking value $1$ on the permutation which sends $i$-th vertex to $j$-th vertex and takes the value zero on other elements of the group. It is a quantum subgroup of $Q^{\rm aut}_{\rm Ban}(\Gamma)$. The surjective CQG morphism $\pi:Q^{\rm aut}_{\rm Ban}(\Gamma)\raro C({\rm Aut}(\Gamma))$ sends the generators to generators. 
% (see definition 3.1 of \cite{Bichon}) An action of a CQG $(\clq,\Delta)$ on a finite graph $\Gamma=(V,E)$ consists of $C^{\ast}$-actions $\alpha:C(V)\raro C(V)\ot \clq$ and $\beta:C(E)\raro C(E)\ot\clq$ such that the following diagram commutes:\\
% \begin{center}\begin{tikzcd}
% C(V)\ot C(V) \arrow{r}{\alpha\underline{\ot}\alpha} \arrow[swap]{d}{m(s_{\ast}\ot t_{\ast})} & C(V)\ot C(V)\ot \clq \arrow{d}{(m(s_{\ast}\ot t_{\ast}))\ot{\rm id}_{\clq}}\\
% C(E) \arrow{r}{\beta} & C(E)\ot \clq
% \end{tikzcd}\end{center}
% where $\alpha\underline{\ot}\alpha$ is the tensor product of the coaction $\alpha$ with itself and $m:C(E)\ot C(E)\raro C(E)$ is the multiplication map. The quantum automorphism group of a finite graph $\Gamma$ is a CQG $(\clq,\Delta)$ acing on $\Gamma$ in the above sense satisfying the following universal property:\\
% \indent If $(A,\Delta_{A})$ is a CQG acting on $\Gamma$ with actions $\alpha^{\prime}:C(V)\raro C(V)\ot A$, $\beta^{\prime}:C(E)\raro C(E)\ot A$, then there is a surjective CQG morphism $\pi:\clq\raro A$ such that $({\rm id}\ot\pi)\alpha=\alpha^{\prime}$ and $({\rm id}\ot\pi)\beta=\beta^{\prime}$
% \edfn

\bthm
(Lemma 3.1.1 of \cite{Fulton}) The quantum automorphism group $(Q^{\rm aut}_{\rm Ban}(\Gamma),\Delta)$ of a finite graph $\Gamma$ with $n$ edges and $m$ vertices (without multiple edge) is the universal $C^{\ast}$-algebra generated by $(u_{ij})_{i,j=1,\ldots,m}$ satisfying the following relations:
\begin{eqnarray}
&& u_{ij}^{\ast}=u_{ij}, u_{ij}u_{ik}=\delta_{jk}u_{ij}, u_{ji}u_{ki}=\delta_{jk}u_{ji}, \sum_{l=1}^{m}u_{il}=\sum_{l=1}^{m}u_{li}=1, 1\leq i,j,k\leq m \label{1}\\
&& u_{s(e_{j})i}u_{t(e_{j})k}=u_{t(e_{j})k}u_{s(e_{j})i}=u_{is(e_{j})}u_{kt(e_{j})}=u_{kt(e_{j})}u_{is(e_{j})}=0, e_{j}\in E, (i,k)\not\in E \label{2}
%&& u_{s(e_{j})s(e_{l})}u_{t(e_{j})t(e_{l})}=u_{t(e_{j})t(e_{l})}u_{s(e_{j})s(e_{l})} \label{3}\\
%&& \sum_{l=1}^{n}u_{s(e_{j})s(e_{l})}u_{t(e_{j})t(e_{l})}=\sum_{l=1}^{n}u_{s(e_{l})s(e_{j})}u_{t(e_{l})t(e_{j})}, 1\leq j\leq n, \label{4}
\end{eqnarray}
where the coproduct on the generators is given by $\Delta(u_{ij})=\sum_{k=1}^{m}u_{ik}\ot u_{kj}$. The $C^{\ast}$-action on the graph is given by $\alpha(p_{i})=\sum_{j}p_{j}\ot u_{ji}$, where $p_{i}$ is the function which takes value $1$ on $i$-the vertex and zero elsewhere.
\ethm
\brmrk
\label{transpose}
Since  $(Q^{\rm aut}_{\rm Ban}(\Gamma),\Delta)$ is a quantum subgroup of $S_{n}^{+}$, it is a Kac algebra and hence $\kappa(u_{ij})=u_{ji}^{\ast}=u_{ji}$. Applying $\kappa$ to the equation $AD=DA$, we get $A^{T}D=DA^{T}$.
\ermrk

With analogy of vertex transitive action of the automorphism group of a graph, we have the following
\bdfn
A graph $\Gamma$ is said to be quantum vertex transitive if the generators $u_{ij}$ of $Q^{\rm aut}_{\rm Ban}(\Gamma)$ are all non zero.
\edfn
\brmrk 
\label{qvertex}
It is easy to see that if a graph is vertex transitive, it must be quantum vertex transitive.
\ermrk
% it must be quantum vertex transitive since the CQG morphism realizing $C({\rm Aut}(\Gamma))$ as a quantum subgroup of $Q^{\rm aut}_{\rm Ban}(\Gamma)$ sends generators to generators and hence any generator $u_{ij}$  of $Q^{\rm aut}_{\rm Ban}(\Gamma)$, which is zero, will be mapped to corresponding $u_{ij}$ of $C({\rm Aut}(\Gamma))$ which will be zero, contradicting the vertex transitivity.
%\ermrk
\bppsn (Corollary 3.7 of \cite{nonlocal1})
\label{ergodic}
The action of $Q_{\rm Ban}^{\rm aut}(\Gamma)$ on $C(V)$ is ergodic if and only if the action is quantum vertex transitive.
\eppsn
\brmrk
For a graph $\Gamma=(V,E)$, when we talk about ergodic action, we always take the corresponding $C^{\ast}$-algebra to be $C(V)$.
\ermrk 
In the next proposition we shall see that in fact for a finite, connected graph $\Gamma$ without multiple edge the CQG $Q_{\rm Ban}^{\rm aut}(\Gamma)$ has a $C^{\ast}$-action on the infinite diemnsional $C^{\ast}$-algebra $C^{\ast}(\Gamma)$.
\bppsn\label{aut} (see Theorem 4.1 of \cite{Web})
Given a directed graph $\Gamma$ without multiple edge, $Q^{\rm aut}_{\rm Ban}(\Gamma)$ has a $C^{\ast}$-action on $C^{\ast}(\Gamma)$. The action is given by 
\begin{eqnarray*}
	&&\alpha(p_{v_{i}})=\sum_{k=1}^{m}p_{v_{k}}\ot u_{ki}\\
	&& \alpha(S_{j})=\sum_{l=1}^{n}S_{l}\ot u_{s(e_{l})s(e_{j})}u_{t(e_{l})t(e_{j})}. 
	\end{eqnarray*}
\eppsn
\bppsn
\label{stateprojection}
Suppose $\Gamma=(V=(v_{1},\ldots,v_{m}),E=(e_{1},\ldots,e_{n}))$ is a finite, directed graph without multiple edges as before. For a ${\rm KMS}_{\beta}$ state $\tau$ on the graph $C^{\ast}$-algebra $C^{\ast}(\Gamma)$, $Q_{\rm Ban}^{\rm aut}(\Gamma)$ preserves $\tau$ if and only if $(\tau\ot{\rm id})\circ\alpha(p_{v_{i}})=\tau(p_{v_{i}})1$ for all $i=1,\ldots,m$.
\eppsn
We start with proving the following Lemma which will be used to prove the Proposition. In the following Lemma, $\Gamma=(V,E)$ is again a finite, directed graph without multiple edges.
\blmma
\label{vertextransitive}
Given any linear functional $\tau$ on $C(V)$, if $Q^{\rm aut}_{\rm Ban}(\Gamma)$ preserves $\tau$, then $\tau(p_{v_{i}})\neq\tau(p_{v_{j}})\Rightarrow u_{ij}=0$.
\elmma
{\it Proof}:\\
Let $i,j$ be such that $\tau(p_{v_{i}})\neq \tau(p_{v_{j}})$. By the assumption,
\begin{eqnarray*}
	&&(\tau\ot{\rm id})\circ\alpha(p_{v_{i}})=\tau(p_{v_{i}})1\\
&\Rightarrow& \sum_{k}\tau(p_{v_{k}})u_{ki}=\tau(p_{v_{i}})1.
\end{eqnarray*}
Multiplying both sides of the last equation by $u_{ji}$ and using the orthogonality, we get $\tau(p_{v_{j}})u_{ji}=\tau(p_{v_{i}})u_{ji}$ i.e. $(\tau(p_{v_{j}})-\tau(p_{v_{i}}))u_{ji}=0$ and hence $u_{ji}=0$ as $\tau(p_{v_{i}})\neq \tau(p_{v_{j}})$. Applying $\kappa$, we get $u_{ij}=u_{ji}=0$.\qed\vspace{0.2in}\\
{\it Proof of  Proposition \ref{stateprojection}}:\\
If $Q^{\rm aut}_{\rm Ban}(\Gamma)$ preserves $\tau$, then $(\tau\ot{\rm id})\circ\alpha(p_{v_{i}})=\tau(p_{v_{i}})1$ for all $i=1,\ldots,m$ trivially. For the converse, given $(\tau\ot{\rm id})\circ\alpha(p_{v_{i}})=\tau(p_{v_{i}})1$ for all $i=1,\ldots,m$, we need to show that $(\tau\ot{\rm id})\circ\alpha(S_{\mu}S_{\nu}^{\ast})=\tau(S_{\mu}S_{\nu}^{\ast})$ for all $\mu,\nu\in E^{\ast}$. The proof is similar to that of Theorem 3.5 of \cite{sou_arn2}.  It is easy to see that for $|\mu|\neq|\nu|$, $(\tau\ot{\rm id})\circ\alpha(S_{\mu}S_{\nu}^{\ast})=0=\tau(S_{\mu}S_{\nu}^{\ast})$. So let $|\mu|=|\nu|$. For $\mu=\nu=e_{i_{1}}e_{i_{2}}\ldots e_{i_{p}}$, we have $S_{\mu}S_{\mu}^{\ast}=S_{i_{1}}\ldots S_{i_{p}}S_{i_{p}}^{\ast}\ldots S_{i_{1}}^{\ast}$. So 
\begin{eqnarray*}
	(\tau\ot{\rm id})\circ\alpha(S_{\mu}S_{\mu}^{\ast})&=&(\tau\ot{\rm id})(\sum S_{j_{1}}\ldots S_{j_{p}}S_{j_{p}}^{\ast}\ldots S_{j_{1}}^{\ast}\ot
	 u_{s(e_{j_{1}})s(e_{i_{1}})}u_{t(e_{j_{1}})t(e_{i_{1}})}\\
	&&\ldots u_{s(e_{j_{p}})s(e_{i_{p}})}u_{t(e_{j_{p}})t(e_{i_{p}})}u_{s(e_{j_{p}})s(e_{i_{p}})}
	\ldots u_{t(e_{j_{1}})t(e_{i_{1}})}u_{s(e_{j_{1}})s(e_{i_{1}})})
	\end{eqnarray*}
By the same argument as given in the proof of the Theorem 3.5 of \cite{sou_arn2}, for $S_{j_{1}}\ldots S_{j_{p}}=0$, $u_{s(e_{j_{1}})s(e_{i_{1}})}u_{t(e_{j_{1}})t(e_{i_{1}})}\ldots u_{s(e_{j_{p}})s(e_{i_{p}})}u_{t(e_{j_{p}})s(t_{i_{p}})}=0$ and hence the last expression equals to
\begin{eqnarray*}
	&& \sum e^{-\beta|\mu|}\tau(p_{t(e_{j_{p}})}) u_{s(e_{j_{1}})s(e_{i_{1}})}u_{t(e_{j_{1}})t(e_{i_{1}})}\ldots u_{s(e_{j_{p}})s(e_{i_{p}})}u_{t(e_{j_{p}})t(e_{i_{p}})}\\
	&&u_{s(e_{j_{p}})s(e_{i_{p}})}\ldots u_{t(e_{j_{1}})t(e_{i_{1}})}u_{s(e_{j_{1}})s(e_{i_{1}})}
\end{eqnarray*}
Observe that any ${\rm KMS}_{\beta}$ state restricts to a state on $C(V)$ so that by Lemma \ref{vertextransitive}, for $\tau(p_{t(e_{j_{p}})})\neq\tau(p_{t(e_{i_{p}})})$, $u_{t(e_{j_{p}})t(e_{i_{p}})}=0$. Using this,  the last summation reduces to 
\begin{eqnarray*}
&& e^{-\beta|\mu|}\sum\tau(p_{t(e_{i_{p}})})u_{s(e_{j_{1}})s(e_{i_{1}})}u_{t(e_{j_{1}})t(e_{i_{1}})}\ldots u_{s(e_{j_{p}})s(e_{i_{p}})}u_{t(e_{j_{p}})t(e_{i_{p}})}\\
 && u_{s(e_{j_{p}})s(e_{i_{p}})}\ldots u_{t(e_{j_{1}})t(e_{i_{1}})}u_{s(e_{j_{1}})s(e_{i_{1}})}.
\end{eqnarray*}
Using  the same arguments used in the proof of Theorem 3.13 in \cite{sou_arn} repeatedly, it can be shown that the last summation actually equals to $e^{-\beta|\mu|}\tau(p_{t(e_{i_{p}})})=\tau(S_{\mu}S_{\mu}^{\ast})$. Hence  
\begin{displaymath}
(\tau\ot{\rm id})\circ\alpha(S_{\mu}S_{\mu}^{\ast})=\tau(S_{\mu}S_{\mu}^{\ast})1.
\end{displaymath}
With similar reasoning it can easily be verified that for $\mu\neq\nu$, $(\tau\ot{\rm id})\circ\alpha(S_{\mu}S_{\nu}^{\ast})=0=\tau(S_{\mu}S_{\nu}^{\ast})1$. Hence by linearity and continuity of $\tau$, for any $a\in C^{\ast}(\Gamma)$, $(\tau\ot{\rm id})\circ\alpha(a)=\tau(a).1$.\qed\vspace{0.1in}\\
If we apply Proposition \ref{stateprojection} to the action of classical automorphism group of a graph on the corresponding graph $C^{\ast}$-algebra, we get the following
\blmma
\label{inv_1}
Given a ${\rm KMS}_{\beta}$ state $\tau$ on $C^{\ast}(\Gamma)$,  we denote the vector $(\tau(p_{v_{1}}),\ldots,\tau(p_{v_{m}}))$ by $\cln^{\tau}$.  If we denote the permutation matrix corresponding to an element $g\in{\rm Aut}(\Gamma)$ by $B$, then ${\rm Aut}(\Gamma)$ preserves $\tau$ if and only if $B\cln^{\tau}=\cln^{\tau}$.
\elmma
{\it Proof}:\\
Follows from the easy observation that $(\tau\ot{\rm id})\circ\alpha(p_{v_{i}})=\tau(p_{v_{i}})1$ implies $B\cln^{\tau}=\cln^{\tau}$ for the classical automorphism group of the graph.\qed.
%\blmma
%\label{gaugecommute}
%For a directed graph $\Gamma$ without multiple edge, the action of the automorphism group commutes with the action of the circle on $C^{\ast}(\Gamma)$.
%\elmma
%{\it Proof}:\\
%$\alpha_{g}(\gamma_{z}(S_{i}))=zS_{j}=\gamma_{z}\alpha_{g}(S_{i})$. By homomorphic property and continuity of both $\alpha_{g}$ and $\gamma_{z}$, $\alpha_{g}\circ\gamma_{z}=\gamma_{z}\circ\alpha_{g}$ for all $z\in \mathbb{T}$ and all $g\in {\rm Aut}(\Gamma)$.\qed
%\blmma
%\label{Weber}
%Let $\Gamma=(E,V)$ be a finite graph without multiple edges or loops and $e_{j}\in E$. Then for any $v\in V$ such that $s^{-1}(v)=\emptyset$, $u_{vs(e_{j})}=u_{s(e_{j})v}=0$ for all $j=1,\ldots,n$, where $u_{ij}$'s are generators of $(Q^{aut}_{Ban}(\Gamma),\Delta)$. 
%\elmma 

\section{Invariance of KMS states under the symmetry of graphs}
In this section we shall study invariance of KMS states of certain classes of graph $C^{\ast}$-algebras under classical and quantum symmetry as mentioned in the introduction. 
%Firstly we shall show that for strongly connected graphs, the unique KMS state at critical inverse temperature is invariant under the action of $Q^{\rm aut}_{\rm Ban}(\Gamma)$. Then we shall look at circulant graphs where the KMS state is not unique, but the classical or quantum automorphism group invariant KMS state is unique. In the last subsection, we shall give an example of a graph $C^{\ast}$-algebra where the KMS state is not unique and though all the KMS states are invariant under the action of classical automorphism group, a unique KMS state is invariant under the action of  $Q^{\rm aut}_{\rm Ban}(\Gamma)$.
%In this section we shall show that imposing the invariance condition on  KMS states under classical or quantum symmetry fixes KMS state for certain class of graph $C^{\ast}$-algebras. More precisely we shall look at three instances. 
\subsection{Strongly connected graphs}
Recall the unique KMS state of $C^{\ast}(\Gamma)$ for a strongly connected graph $\Gamma$ with vertex matrix $D$. We denote the $ij$-th entry of $D$ by $d_{ij}$. The unique KMS state at the inverse critical temperature ${\rm ln}(\rho(D))$ is determined by the unique normalized Perron-Frobenius eigen vector of $D$ corresponding to the eigen value $\rho(D)$. If the state is denoted by $\tau$, the eigen vector is given by $((\tau(p_{v_{i}})))_{i=1,\ldots,m}$ where $m$ is the number of vertices. Now recall the action of $Q^{\rm aut}_{\rm Ban}(\Gamma)$ on $C^{\ast}(\Gamma)$. We continue to denote the matrix $((u_{ij}))$ by $A$.
\blmma
\label{0state}
$(\tau\ot {\rm id})\circ\alpha(p_{v_{i}})=\tau(p_{v_{i}})1 \forall \ i=1,\ldots,m$.
\elmma
{\it Proof}:\\ 
For a state $\phi$ of $Q_{\rm aut}^{\rm Ban}(\Gamma)$, we denote the vector whose $i$-th entry is $(\tau\ot\phi)\circ\alpha(p_{v_{i}})$ by ${\bf v}_{\phi}$. Then
\begin{eqnarray*}
	(D{\bf v}_{\phi})_{i}&=&\sum_{j}d_{ij}(\tau\ot\phi)\circ\alpha(p_{v_{j}})\\
	&=& \sum_{j,k}\tau(p_{v_{k}})\phi(d_{ij}u_{kj})\\
	&=& \sum_{k}\tau(p_{v_{k}})\phi(\sum_{j}d_{ij}u_{kj})\\
	&=& \sum_{k}\tau(p_{v_{k}})\phi(\sum_{j}u_{ji}d_{jk}) \  \ (DA^{T}=A^{T}D  \ by \ remark \ \ref{transpose})\\
	&=& \sum_{j,k}d_{jk}\tau(p_{v_{k}})\phi(u_{ji})\\
	&=&\rho(D)\sum_{j}\tau(p_{j})\phi(u_{ji}) \  \ (D((\tau(p_{i})))^{T}=\rho(D)((\tau(p_{v_{i}})))^{T})\\
	&=& \rho(D) (\tau\ot\phi)\circ\alpha(p_{v_{i}}). 
	\end{eqnarray*}
Hence ${\bf v}_{\phi}$ is an eigen vector of $D$ corresponding to the eigen value $\rho(D)$. By the one dimensionality of the eigen space we have some constant $C_{\phi}$ such that $(\tau\ot \phi)\circ\alpha(p_{v_{i}})=C_{\phi}\tau(p_{v_{i}})$ for all $i=1,\ldots,m$. To determine the constant $C_{\phi}$ we take the summation over $i$ on both sides and get
\begin{eqnarray*}
	&&\sum_{i}(\tau\ot\phi)\circ\alpha(p_{v_{i}})=C_{\phi}\sum_{i}\tau(p_{v_{i}})\\
	&\Rightarrow& \sum_{i,j}\tau(p_{v_{j}})\phi(u_{ji})=C_{\phi}\\
	&\Rightarrow& \sum_{j}\tau(p_{v_{j}})\sum_{i}\phi(u_{ji})=C_{\phi}\\
	&\Rightarrow& C_{\phi}=1.
	\end{eqnarray*}
From above calculation it is clear that $C_{\phi}=1$ for all states $\phi$ on $Q^{\rm aut}_{\rm Ban}(\Gamma)$. Hence $\phi((\tau\ot{\rm id})\circ\alpha(p_{v_{i}}))=\phi(\tau(p_{v_{i}})1)$ for all $i$ and for all state $\phi$ which implies that $(\tau\ot{\rm id})\circ\alpha(p_{v_{i}})=\tau(p_{v_{i}})1$ for all $i=1,\ldots,m$.\qed
\bppsn
\label{statepreserve}
For a strongly connected graph $\Gamma$, $Q^{\rm aut}_{\rm Ban}(\Gamma)$ preserves the unique KMS state of $C^{\ast}(\Gamma)$.
\eppsn
{\it Proof}:\\
This follows from Lemma \ref{0state} and Proposition \ref{stateprojection}.\qed
\brmrk
In light of the above Proposition \ref{statepreserve}, for strongly connected graphs, we can relax the condition on the graph in \cite{sou_arn2}. In that paper it was assumed that all the row sums of the vertex matrix have to be equal to ensure that $Q^{\rm aut}_{\rm Ban}(\Gamma)$ belongs to the category $\clc^{\Gamma}_{\tau}$ (see \cite{sou_arn2} for notation)  for the unique KMS state $\tau$ on a strongly connected graph $\Gamma$. Now we have  for a strongly connected graph $\Gamma$ with its unique KMS state $\tau$, the category $\clc_{\tau}^{\Gamma}$  contains $Q^{\rm aut}_{\rm Ban}(\Gamma)$ and for that one does not require  all the row sums to be equal.
\ermrk
We end this subsection with a proposition about the non ergodicity of the action of $Q_{\rm Ban}^{\rm aut}(\Gamma)$ on $C(V)$ for a strongly connected graph $\Gamma=(V,E)$. Note that the following proposition does not deal with states on the infinite dimensional $C^{\ast}$-algebra $C^{\ast}(\Gamma)$.
\bppsn
\label{erg}
%1.The action of $Q_{\rm aut}^{\rm Ban}(\Gamma)$ on a strongly connected graph $\Gamma$ is not quantum vertex transitive if  the Perron-Frobenius eigen space of the vertex matrix is not spanned by the vector $(1,\ldots,1)$.\\
 For a strongly connected graph $\Gamma$, if the Perron-Frobenius eigen vector is not spanned by $(1,\ldots,1)$, then the action of $Q^{\rm aut}_{\rm Ban}(\Gamma)$ is non ergodic.
\eppsn
{\it Proof}:\\
It follows from Lemma \ref{vertextransitive}, Lemma \ref{0state} and Proposition \ref{ergodic}.\qed
%Statement of 1 follows directly from Corollary \ref{vertextransitive}. Statement 2 follows from statement 1 and Proposition \ref{ergodic}.\qed
\subsection{Circulant graphs}
Consider a finite graph $\Gamma=(V,E)$ with $m$-vertices $(v_{1},\ldots,v_{m})$. Recall the notation $\cln^{\tau}$ for a ${\rm KMS}_{\beta}$ state $\tau$ on $C^{\ast}(\Gamma)$.
\blmma
\label{trans}
Given a graph $\Gamma$ such that ${\rm Aut}(\Gamma)$ acts transitively on its vertices, $B\cln^{\tau}=\cln^{\tau}$ for all $B\in{\rm Aut}(\Gamma)$ if and only if  $\cln^{\tau}_{i}=\cln^{\tau}_{j}$ for all $i,j=1,\ldots,m$.
\elmma
{\it Proof}:\\
If $\cln^{\tau}_{i}=\cln^{\phi}_{j}$ for all $i,j=1,\ldots,m$, then $B\cln^{\tau}=\cln^{\tau}$ for all $B\in{\rm Aut}(\Gamma)$ trivially. For the converse, let $\cln^{\tau}_{i}\neq \cln^{\tau}_{j}$ for some $i,j$. Since the action of the automorphism group is transitive, there is some $B\in{\rm Aut}(\Gamma)$ so that $B(v_{i})=v_{j}$ and hence $B\cln^{\tau}\neq \cln^{\tau}$.\qed\vspace{0.1in}\\
Now recall from the discussion following Corollary \ref{not_unique} and from Lemma \ref{circonetemp} that for a circulant graph, the KMS state exists only at the critical inverse temperature, but it is not necessarily unique. We shall prove that if we further assume the invariance of such a state under the action of the automorphism group of the graph, then it is unique. 
\bppsn
For a circulant graph $\Gamma$ there exists a unique ${\rm Aut}(\Gamma)$ invariant KMS state on $C^{\ast}(\Gamma)$. 
\eppsn
{\it Proof}:\\
 Since  for a circulant graph, the automorphism group acts transitively on the set of vertices, by Lemma \ref{trans} and Lemma \ref{inv_1},   a ${\rm KMS}_{\beta}$ state $\tau$ is ${\rm Aut}(\Gamma)$ invariant if and only if $\tau(p_{v_{i}})=\frac{1}{m}$ for all $i$. This coupled with the fact that $(\frac{1}{m},\ldots,\frac{1}{m})$ is an eigen vector corresponding to the eigen value $\lambda$ ($=$   spectral radius) finishes the proof of the proposition.\qed
\brmrk
\label{future}
The group invariant KMS state is also invariant under the action of quantum automorphism group of the underlying graph. Since $\tau(p_{v_{i}})=\tau(p_{v_{j}})$, it is easy to see that for the action of $Q^{\rm aut}_{\rm Ban}(\Gamma)$, $(\tau\ot{\rm id})\circ\alpha(p_{v_{i}})=\tau(p_{v_{i}})1$ for all $i=1,\ldots,m$. Hence an application of Proposition \ref{stateprojection} finishes the proof of the claim.
\ermrk

\subsection{Graph of Mermin-Peres magic square game}
We start this subsection by clarifying a few notations to be used in this subsection. Given an undirected graph $\Gamma$, we make it directed by declaring that both $(i,j)$  and $(j,i)$ are in the edge set whenever there is an edge between two vertices $v_{i}$ and $v_{j}$. The vertex matrix of such directed graph is symmetric by definition. In this subsection we use the notation $\overrightarrow{\Gamma}$ for the directed graph coming from an undirected graph $\Gamma$ in this way. $\Gamma$ will always denote an undirected graph.
\brmrk
\label{arrow}
By definition, $Q^{\rm aut}_{\rm Ban}(\overrightarrow{\Gamma})\cong Q^{\rm aut}_{\rm Ban}(\Gamma)$ and hence ${\rm Aut}(\overrightarrow{\Gamma})\cong{\rm Aut}(\Gamma)$ (see \cite{Ban}).
\ermrk Given two  graphs $\Gamma_{1}=(V_{1},E_{1}),\Gamma_{2}=(V_{2},E_{2})$  their disjoint union $\Gamma_{1}\cup\Gamma_{2}$ is defined to be the graph $\Gamma=(V,E)$ such that $V=V_{1}\cup V_{2}$. There is an edge between two vertices $v_{i}, v_{j}\in V_{1}\cup V_{2}$ if both the vertices belong to either $\Gamma_{1}$ or $\Gamma_{2}$ and they have an edge in the corresponding graph. 
\bppsn
\label{iso}
Let $\Gamma_{1}, \Gamma_{2}$ be two non isomorphic  connected graphs. Then the automorphism group of $\overrightarrow{\Gamma_{1}\cup\Gamma_{2}}$ is given by ${\rm Aut}(\overrightarrow{\Gamma_{1}})\times{\rm Aut}(\overrightarrow{\Gamma_{2}})$.
\eppsn
{\it Proof}:\\
Combining Theorem 2.5 of \cite{thesis} and Remark \ref{arrow}, the result follows. \qed

\bppsn
\label{infKMS}
Let $\Gamma_{1}$ and $\Gamma_{2}$ be two non isomorphic  connected graphs such that $\overrightarrow{\Gamma_{1}}$ and $\overrightarrow{\Gamma_{2}}$ has   symmetric vertex matrices $D_{1}$ and $D_{2}$ having equal spectral radius say $\lambda$. Then for the graph $\Gamma=\Gamma_{1}\cup \Gamma_{2}$, $C^{\ast}(\overrightarrow{\Gamma})$ has infinitely many KMS states at the inverse critical temperature ${\rm ln}(\lambda)$ such that all of them are invariant under the action of ${\rm Aut}(\overrightarrow{\Gamma})\cong{\rm Aut}(\overrightarrow{\Gamma_{1}})\times {\rm Aut}(\overrightarrow{\Gamma_{2}})$.
\eppsn
To prove the proposition we require the following 
\blmma
\label{joint}
Let $A\in M_{n}(\mathbb{C})$ and $B\in M_{m}(\mathbb{C})$. Then the spectral radius of the matrix $\begin{bmatrix}
A & 0_{n\times m}\\
0_{m\times n} & B
\end{bmatrix}$ is equal to ${\rm max}\{{\rm sp}(A),{\rm sp}(B)\}$.
\elmma
{\it Proof}:\\
It follows from the simple observation that any eigen value of the matrix  $\begin{bmatrix}
A & 0_{n\times m}\\
0_{m\times n} & B
\end{bmatrix}$ is either an eigen value of $A$ or an eigen value of $B$.\qed\vspace{0.2in}\\
{\it Proof of Proposition \ref{infKMS}}:\\
We assume that $\overrightarrow{\Gamma_{1}}$ has $n$-vertices and $\overrightarrow{\Gamma_{2}}$ has $m$-vertices. Let us denote the vertex matrix of $\overrightarrow{\Gamma}$ by $D$. $D$ is given by the matrix \begin{center}$\begin{bmatrix}
D_{1} & 0_{n\times m}\\
0_{m\times n} & D_{2}
\end{bmatrix}$.\end{center} Then the spectral radius is equal to $\lambda$ by Lemma \ref{joint}. Also it is easy to see that the spectral radius is an eigen value of the matrix $D$. Now $\lambda$ has one dimensional eigen space for $D_{1}$ spanned by say $w_{1}$ and  one dimensional eigen space for $D_{2}$ spanned by say $w_{2}$ as both the graphs are connected and hence strongly connected as directed graphs. We take both the eigen vectors normalized for convenience. Then for $D$, the eigen space corresponding to $\lambda$ is two dimensional spanned by the vectors ${\bf w_{1}}=(w_{1},0_{m})$ and ${\bf w_{2}}=(0_{n},w_{2})$ where $0_{k}$ is the zero $k$-tuple. So the eigen space of $D$ corresponding to the eigen value $\lambda$ is given by $\{\xi{\bf w_{1}}+\eta{\bf w_{2}}:(\xi,\eta)\in\mathbb{C}^2- ( 0,0) \}$. For any $(\xi,\eta)\in\mathbb{C}^{2}$, $\xi{\bf w_{1}}+\eta{\bf w_{2}}=(\xi w_{1},\eta w_{2})$.  It is easy to see that there are infinitely many choice of $\xi,\eta$ such that corresponding eigen vector is normalized with all its entries non negative  which in turn give rise to infinitely many  KMS states. The set of normalized vectors is given by $\{(\xi w_{1},(1-\xi)w_{2}):0\leq\xi\leq 1\}$. We shall show that any $B\in{\rm Aut}(\overrightarrow{\Gamma})$ keeps such a normalized eigen vector invariant. Let ${\bf w}=(\xi w_{1},(1-\xi)w_{2})$ be one such choice. By Proposition \ref{iso},  any $B\in {\rm Aut}(\overrightarrow{\Gamma})$ can be written in the matrix form $\begin{bmatrix}
B_{1} & 0_{n\times m}\\
0_{m\times n} & B_{2}.\end{bmatrix}$, for $B_{i}\in {\rm Aut}(\overrightarrow{\Gamma_{i}})$ and $i=1,2$. Then $B{\bf w}=(\xi B_{1}w_{1},(1-\xi) B_{2}w_{2})$. Since $w_{1}, w_{2}$ are Perron-Frobenius eigen vectors of $D_{1}, D_{2}$  respectively with  both the graphs $\overrightarrow{\Gamma_{1}}, \overrightarrow{\Gamma_{2}}$ strongly connected,  by Lemma \ref{0state}, $B_{i}(w_{i})=w_{i}$ for $i=1,2$. So  $B{\bf w}=(\xi w_{1},(1-\xi)w_{2})={\bf w}$. Now an application of  Lemma \ref{inv_1} completes the proof of the proposition.\qed\vspace{0.1in}\\
Now we turn to the main object of study of this Subsection.
%\bppsn
%\label{infKMS}
%For two strongly connected non isomorphic graphs $\Gamma_{1}$ and $\Gamma_{2}$ with equal spectral radius $\lambda$, their disjoint union $\Gamma_{1}\cup\Gamma_{2}$ has infinitely many KMS states occurring at the critical inverse temperature ${\rm ln}(\lambda)$. All the KMS states are invariant under the action of ${\rm Aut}(\Gamma_{1}\cup\Gamma_{2})\cong {\rm Aut}(\Gamma_{1})\times{\rm Aut}(\Gamma_{2})$.
%\eppsn 

\subsubsection{Linear Binary Constraint System (LBCS)}
A linear binary constraint system (LBCS) $\clf$ consists of a family of binary variables $x_{1},\ldots,x_{n}$ and constraints $C_{1},\ldots,C_{m}$, where each $C_{l}$ is a linear equation over $\mathbb{F}_{2}$ in some subset of the variables i.e. each $C_{l}$ is of the form $\sum_{x_{i}\in S_{l}}x_{i}=b_{l}$ for some $S_{l}\subset\{x_{1},\ldots,x_{n}\}$. An LBCS is said to be {\it satisfiable} if there is an assignment of values from $\mathbb{F}_{2}$ to the variables $x_{i}$ such that every constraint $C_{l}$ is satisfied. For every LBCS there is an associated LBCS game. For details on LBCS games readers are referred to \cite{nonlocal1} and \cite{nonlocal2}. For the following definition we need the concept of perfect quantum strategy for nonlocal games. We shall not discuss it here and readers are again referred to \cite{nonlocal1}, \cite{nonlocal2} for details on quantum strategy.
\bdfn
An LBCS is said to be quantum satisfiable if there exists a perfect quantum strategy for the corresponding LBCS game.
\edfn
Now we shall give an example of an LBCS which is quantum satisfiable but not satisfiable. Consider the following LBCS:
\begin{eqnarray*}
	&& x_{1}+x_{2}+x_{3}=0 \ \ \ \ \ \ \ \ \ \ \ x_{1}+x_{4}+x_{7}=0\\
	&& x_{4}+x_{5}+x_{6}=0 \ \ \ \ \ \ \ \ \ \ \ x_{2}+x_{5}+x_{8}=0\\
	&& x_{7}+x_{8}+x_{9}=0 \ \ \ \ \ \ \ \ \ \ \ x_{3}+x_{6}+x_{9}=1,
\end{eqnarray*} 
where the addition is over $\mathbb{F}_{2}$. It is easy to see that the above LBCS is not satisfiable since summing up all equations modulo $2$ we get $0=1$. The fact that it is quantum satisfiable is more non trivial. For that we refer the reader to \cite{nonlocal2}. The game corresponding to the above LBCS is called the Mermin-Peres magic square game. Corresponding to every LBCS $\clf$ one can associate a graph (see section 6.2 of \cite{nonlocal2}) to be denoted by $\clg(\clf)$. 

\bdfn
\label{homo}
Given an LBCS $\clf$, its homogenization $\clf_{0}$ is defined to be the LBCS obtained by assigning zero to the right hand side of every constraint $C_{l}$.
\edfn
 We have the following (Theorem 6.2 and 6.3 of \cite{nonlocal2})
\bthm
\label{quantum}
Given an LBCS $\clf$, $\clf$ is (quantum) satisfiable if and only if the graphs $\clg(\clf)$ and $\clg(\clf_{0})$ are (quantum) isomorphic.
\ethm
 From now on $\clf$ will always mean the LBCS corresponding to Mermin-Peres magic square game. In light of the  Theorem \ref{quantum} and the discussion just before the Definition \ref{homo} we readily see that for the Mermin-Peres magic square game, the corresponding graphs $\clg(\clf)$ and $\clg(\clf_{0})$ are quantum isomorphic, but not isomorphic.  Both the graphs are vertex transitive (in fact they are Cayley as mentioned in \cite{nonlocal1}) and hence quantum vertex transitive by Remark \ref{qvertex}. Combining the facts that $\clg(\clf)$ and $\clg(\clf_{0})$ are quantum isomorphic and quantum vertex transitive with Lemma 4.11 of \cite{nonlocal1}, we get 
\blmma
\label{qvt}
For the LBCS $\clf$ corresponding to the Mermin Peres magic square game, the disjoint union of $\clg(\clf)$ and $\clg(\clf_{0})$ is quantum vertex transitive.
\elmma
By Remark \ref{arrow},
\bcrlre
\label{qvt1}
$\overrightarrow{\clg(\clf)\cup\clg(\clf_{0})}$ is quantum vertex transitive.
\ecrlre

\begin{tikzpicture}[scale=0.5, transform shape]
\tikzmath{\x1 = 0; \y1 =2; \y2=-10;}

\path (0,2) node[circle,draw](1){000};
\path (2,2) node[circle,draw](2){011};
\path (4,2) node[circle,draw](3){101};
\path (6,2) node[circle,draw](4){110};

\path (10,2) node[circle,draw](5){000};
\path (12,2) node[circle,draw](6){011};
\path (14,2) node[circle,draw](7){101};
\path (16,2) node[circle,draw](8){110};

\path (20,2) node[circle,draw](9){000};
\path (22,2) node[circle,draw](10){011};
\path (24,2) node[circle,draw](11){101};
\path (26,2) node[circle,draw](12){110};
%%%%%%%
\path (0,-10) node[circle,draw](13){000};
\path (2,-10) node[circle,draw](14){011};
\path (4,-10) node[circle,draw](15){101};
\path (6,-10) node[circle,draw](16){110};

\path (10,-10) node[circle,draw](17){000};
\path (12,-10) node[circle,draw](18){011};
\path (14,-10) node[circle,draw](19){101};
\path (16,-10) node[circle,draw](20){110};

\path (20,-10) node[circle,draw](21){111};
\path (22,-10) node[circle,draw](22){100};
\path (24,-10) node[circle,draw](23){010};
\path (26,-10) node[circle,draw](24){001};
%%%%%%%
\draw (1) .. controls (1,3) and (3,3) .. (3);
\draw (1) .. controls (1,3.5) and (5,3.5) .. (4);
\draw (2) .. controls (3,3) and (5,3) .. (4);

\draw (5) .. controls (11,3) and (13,3) .. (7);
\draw (5) .. controls (11,3.5) and (15,3.5) .. (8);
\draw (6) .. controls (13,3) and (15,3) .. (8);

\draw (9) .. controls (21,3) and (23,3) .. (11);
\draw (9) .. controls (21,3.5) and (25,3.5) .. (12);
\draw (10) .. controls (23,3) and (25,3) .. (12);
%%%%%%%
\draw (13) .. controls (1,-11) and (3,-11) .. (15);
\draw (13) .. controls (1,-11.5) and (5,-11.5) .. (16);
\draw (14) .. controls (3,-11) and (5,-11) .. (16);

\draw (17) .. controls (11,-11) and (13,-11) .. (19);
\draw (17) .. controls (11,-11.5) and (15,-11.5) .. (20);
\draw (18) .. controls (13,-11) and (15,-11) .. (20);

\draw (21) .. controls (21,-11) and (23,-11) .. (23);
\draw (21) .. controls (21,-11.5) and (25,-11.5) .. (24);
\draw (22) .. controls (23,-11) and (25,-11) .. (24);

%%%%%%%
\draw (1) -- (2);
\draw (2) -- (3);
\draw (3) -- (4);

\draw (5) -- (6);
\draw (6) -- (7);
\draw (7) -- (8);

\draw (9) -- (10);
\draw (10) -- (11);
\draw (11) -- (12);

\draw (13) -- (14);
\draw (14) -- (15);
\draw (15) -- (16);

\draw (17) -- (18);
\draw (18) -- (19);
\draw (19) -- (20);

\draw (21) -- (22);
\draw (22) -- (23);
\draw (23) -- (24);
%%%%%%%%
\foreach \x in {15,16,19,20,21,22}
\draw (1) -- (\x);
\foreach \x in {15,16,17,18,23,24}
\draw (2) -- (\x);
\foreach \x in {13,14,19,20,23,24}
\draw (3) -- (\x);
\foreach \x in {13,14,17,18,21,22}
\draw (4) -- (\x);

\foreach \x in {14,16,18,20,21,23}
\draw (5) -- (\x);
\foreach \x in {14,16,17,19,22,24}
\draw (6) -- (\x);
\foreach \x in {13,15,18,20,22,24}
\draw (7) -- (\x);
\foreach \x in {13,15,17,19,21,23}
\draw (8) -- (\x);

\foreach \x in {14,15,18,19,21,24}
\draw (9) -- (\x);
\foreach \x in {14,16,17,20,22,23}
\draw (10) -- (\x);
\foreach \x in {13,16,18,19,22,23}
\draw (11) -- (\x);
\foreach \x in {13,16,17,20,21,24}
\draw (12) -- (\x);
%%%%%%
\node at (3,3.8) {\Large{$x_{1}+x_2+x_3=0$}};
\node at (13,3.8) {\Large{$x_{4}+x_5+x_6=0$}};
\node at (23,3.8) {\Large{$x_{7}+x_8+x_9=0$}};
%%%%%%%%
\node at (3,-11.8) {\Large{$x_{1}+x_4+x_7=0$}};
\node at (13,-11.8) {\Large{$x_{2}+x_5+x_8=0$}};
\node at (23,-11.8) {\Large{$x_{3}+x_6+x_9=1$}};
\end{tikzpicture}
\begin{center}
	{Figure 1: Graph $\clg(\clf)$}
	\end{center}
\begin{tikzpicture}[scale=0.5, transform shape]
\tikzmath{\x1 = 0; \y1 =2; \y2=-10;}

\path (0,2) node[circle,draw](1){000};
\path (2,2) node[circle,draw](2){011};
\path (4,2) node[circle,draw](3){101};
\path (6,2) node[circle,draw](4){110};

\path (10,2) node[circle,draw](5){000};
\path (12,2) node[circle,draw](6){011};
\path (14,2) node[circle,draw](7){101};
\path (16,2) node[circle,draw](8){110};

\path (20,2) node[circle,draw](9){000};
\path (22,2) node[circle,draw](10){011};
\path (24,2) node[circle,draw](11){101};
\path (26,2) node[circle,draw](12){110};
%%%%%%%
\path (0,-10) node[circle,draw](13){000};
\path (2,-10) node[circle,draw](14){011};
\path (4,-10) node[circle,draw](15){101};
\path (6,-10) node[circle,draw](16){110};

\path (10,-10) node[circle,draw](17){000};
\path (12,-10) node[circle,draw](18){011};
\path (14,-10) node[circle,draw](19){101};
\path (16,-10) node[circle,draw](20){110};

\path (20,-10) node[circle,draw](21){111};
\path (22,-10) node[circle,draw](22){100};
\path (24,-10) node[circle,draw](23){010};
\path (26,-10) node[circle,draw](24){001};
%%%%%%%
\draw (1) .. controls (1,3) and (3,3) .. (3);
\draw (1) .. controls (1,3.5) and (5,3.5) .. (4);
\draw (2) .. controls (3,3) and (5,3) .. (4);

\draw (5) .. controls (11,3) and (13,3) .. (7);
\draw (5) .. controls (11,3.5) and (15,3.5) .. (8);
\draw (6) .. controls (13,3) and (15,3) .. (8);

\draw (9) .. controls (21,3) and (23,3) .. (11);
\draw (9) .. controls (21,3.5) and (25,3.5) .. (12);
\draw (10) .. controls (23,3) and (25,3) .. (12);
%%%%%%%
\draw (13) .. controls (1,-11) and (3,-11) .. (15);
\draw (13) .. controls (1,-11.5) and (5,-11.5) .. (16);
\draw (14) .. controls (3,-11) and (5,-11) .. (16);

\draw (17) .. controls (11,-11) and (13,-11) .. (19);
\draw (17) .. controls (11,-11.5) and (15,-11.5) .. (20);
\draw (18) .. controls (13,-11) and (15,-11) .. (20);

\draw (21) .. controls (21,-11) and (23,-11) .. (23);
\draw (21) .. controls (21,-11.5) and (25,-11.5) .. (24);
\draw (22) .. controls (23,-11) and (25,-11) .. (24);
%%%%%%%
\draw (1) -- (2);
\draw (2) -- (3);
\draw (3) -- (4);

\draw (5) -- (6);
\draw (6) -- (7);
\draw (7) -- (8);

\draw (9) -- (10);
\draw (10) -- (11);
\draw (11) -- (12);

\draw (13) -- (14);
\draw (14) -- (15);
\draw (15) -- (16);

\draw (17) -- (18);
\draw (18) -- (19);
\draw (19) -- (20);

\draw (21) -- (22);
\draw (22) -- (23);
\draw (23) -- (24);
%%%%%%%%
\foreach \x in {15,16,19,20,23,24}
\draw (1) -- (\x);
\foreach \x in {15,16,17,18,20,21}
\draw (2) -- (\x);
\foreach \x in {13,14,19,20,21,22}
\draw (3) -- (\x);
\foreach \x in {13,14,17,18,23,24}
\draw (4) -- (\x);

\foreach \x in {14,16,18,20,22,24}
\draw (5) -- (\x);
\foreach \x in {14,16,17,19,21,23}
\draw (6) -- (\x);
\foreach \x in {13,15,18,20,21,23}
\draw (7) -- (\x);
\foreach \x in {13,15,17,19,22,24}
\draw (8) -- (\x);

\foreach \x in {14,15,18,19,22,23}
\draw (9) -- (\x);
\foreach \x in {14,15,17,20,21,24}
\draw (10) -- (\x);
\foreach \x in {13,16,18,19,21,24}
\draw (11) -- (\x);
\foreach \x in {13,16,17,20,22,23}
\draw (12) -- (\x);
%%%%%%
\node at (3,3.8) {\Large{$x_{1}+x_2+x_3=0$}};
\node at (13,3.8) {\Large{$x_{4}+x_5+x_6=0$}};
\node at (23,3.8) {\Large{$x_{7}+x_8+x_9=0$}};
%%%%%%%%
\node at (3,-11.8) {\Large{$x_{1}+x_4+x_7=0$}};
\node at (13,-11.8) {\Large{$x_{2}+x_5+x_8=0$}};
\node at (23,-11.8) {\Large{$x_{3}+x_6+x_9=0$}};
\end{tikzpicture}
\begin{center}
	{Figure 2: Graph $\clg(\clf_{0})$}
\end{center}
It can be verified that both the graphs $\clg(\clf)$ and $\clg(\clf_{0})$ are  connected with $24$ vertices each such that the vertex matrices of the graphs $\overrightarrow{\clg(\clf)}$ and $\overrightarrow{\clg(\clf)_{0}}$ have spectral radius $9$. Then since they are non isomorphic, by Proposition \ref{infKMS}, $C^{\ast}(\overrightarrow{\clg(\clf)\cup\clg(\clf_{0})})$ has infinitely many KMS states at the inverse critical temperature ${\rm ln}(9)$ all of which are invariant under the action of the classical automorphism group of $\overrightarrow{\clg(\clf)\cup\clg(\clf_{0})}$.  But as mentioned earlier, if we further assume that the KMS state at critical temperature is invariant under the action of $Q^{\rm aut}_{\rm Ban}(\overrightarrow{\clg(\clf)\cup\clg(\clf_{0})})$, then it is necessarily unique. We prove it in the next theorem.
\bthm
For the LBCS $\clf$ corresponding to the Mermin-Peres magic square game, \\
$C^{\ast}(\overrightarrow{\clg(\clf)\cup\clg(\clf_{0})})$ has a unique $ Q^{\rm aut}_{\rm Ban}(\overrightarrow{\clg(\clf)\cup\clg(\clf_{0})})$-invariant KMS state $\tau$ given by
\begin{displaymath}
\tau(S_{\mu}S_{\nu}^{\ast})=\delta_{\mu,\nu}\frac{1}{9^{|\mu|}48}.
\end{displaymath}

\ethm
{\it Proof}:\\
Since the graph $\overrightarrow{(\clg(\clf)\cup\clg(\clf_{0})})$ is quantum vertex transitive by Corollary \ref{qvt1}, for any KMS state $\tau$ on $C^{\ast}(\overrightarrow{\clg(\clf)\cup\clg(\clf_{0})}))$ at the critical inverse temperature ${\rm ln}(9)$ which is preserved by $Q^{\rm aut}_{\rm Ban}(\overrightarrow{\clg(\clf)\cup\clg(\clf_{0})})$, $\tau(p_{v_{i}})=\tau(p_{v_{j}})$ for all $i,j$ (see Lemma \ref{vertextransitive}). That forces $\tau(p_{v_{i}})$ to be $\frac{1}{48}$ for all $i=1,\ldots,48$. Since $(\frac{1}{48},\ldots,\frac{1}{48})$ is an eigen vector corresponding to the eigen value $9$, there is a unique KMS state on $C^{\ast}(\overrightarrow{\clg(\clf)\cup\clg(\clf_{0})})$ at the inverse critical temperature ${\rm ln}(9)$ satisfying (see (a) of \ref{exist_KMS})
\begin{displaymath}
\tau(S_{\mu}S_{\nu}^{\ast})=\delta_{\mu,\nu}\frac{1}{9^{|\mu|}48}.
\end{displaymath}
$Q^{\rm aut}_{\rm Ban}(\overrightarrow{\clg(\clf)\cup\clg(\clf_{0})})$ preserves the above KMS state by following the same line of arguments as given in Remark \ref{future}.  To complete the proof, we need to show that the only possible inverse temperature where a KMS state could occur is the critical inverse temperature. For that  first notice that the graph $\overrightarrow{\clg(\clf)\cup\clg(\clf_{0})}$ is without sink. We have already observed that the spectral radius $9$ has an eigen vector with all entries strictly positive (column vector with all its entries $\frac{1}{48}$). Also the vertex matrix of the graph $\overrightarrow{\clg(\clf)\cup\clg(\clf_{0})}$ is symmetric implying that $(\frac{1}{48},\ldots,\frac{1}{48})D=9(\frac{1}{48},\ldots,\frac{1}{48})$. Hence an application of Lemma \ref{onetemp} finishes the proof of the theorem. \qed
\section*{Concluding remarks}
1. We conjecture that the converse of the Proposition \ref{erg} is true, i.e.  a strongly connected graph $\Gamma$ is quantum vertex transitive if and only if the Perron-Frobenius eigen space is spanned by the vector $(1,\ldots,1)$. Even if the conjecture is false it seems that it is hard to find an example of a strongly connected graph with `small' number of vertices whose Perron-Frobenius eigen vector is spanned by $(1,\ldots,1)$, but the graph fails to be quantum vertex transitive. \vspace{0.1in}\\
2. In all the examples considered in this paper, KMS states always occur at the inverse critical temperature. But in general it might be interesting to see if some natural symmetry could also fix the inverse temperature. In this context, one can possibly look at the graphs with sink which has richer supply of KMS states (see \cite{watatani}).\vspace{0.1in}\\
{\bf Acknowledegement}: The first author acknowledges support from Department of Science and Technology, India (DST/INSPIRE/04/2016/002469). The second author acknowledges support from Science and Engineering Research Board, India (PDF/2017/001795 ). Both the authors would like to thank Malay Ranjan Biswal, Sruthi C. K. for helping us to draw the figures as well as to find eigen values of the vertex matrices coming from the Mermin-Peres magic square game using Python.

Soumalya Joardar \\
Theoretical Science Unit,\\ 
JNCASR, Bangalore-560064, India\\ 
email: soumalya.j@gmail.com 
\vspace{0.1in}\\
Arnab Mandal\\
School Of Mathematical Sciences\\
NISER, HBNI,  Bhubaneswar,  Jatni-752050, India\\
email: arnabmaths@gmail.com
\end{document}